\documentclass{amsart}
\usepackage{amssymb}
\usepackage{amsmath}
\usepackage{amsfonts}
\usepackage{graphicx}
\newtheorem{theorem}{Theorem}[section]

\newtheorem{prop}[theorem]{Proposition}

\newtheorem{remark}[theorem]{Remark}

\begin{document}

\title{Nonconcentration of energy for a semilinear Skyrme model}

\author{Dan-Andrei Geba, S. G. Rajeev} 

\address{Department of Mathematics,
University of Rochester, Rochester, NY 14627, USA}
\email{dangeba@math.rochester.edu}

\address{Department of Physics and Astronomy,  Department of Mathematics,
University of Rochester,
Rochester, NY 14627, USA}
\email{rajeev@pas.rochester.edu}

\subjclass[2000]{35L71, 81T13}
\keywords{wave maps; Skyrme model; global solutions}

\date{}

\begin{abstract}
We continue our investigation \cite{GR1} of a model introduced by Adkins and Nappi \cite{AN}, in which $\omega$ mesons stabilize chiral solitons. The aim of this article is to show that the energy associated to equivariant solutions does not concentrate.
\end{abstract}

\maketitle

\section{Introduction}

Physicists worry a great deal about the divergences of quantum field theories: realistic theories make accurate predictions only after a mysterious process of renormalization that removes divergences. Much of the motivation for new theories of elementary particles and gravity is the search for more fundamental theories which are finite.

It is perhaps less appreciated among physicists that classical physics often leads to divergences (or singularities) as well. The laws of classical physics are encoded as partial differential equations (PDEs). They do not always have regular solutions, even when the boundary (or initial) conditions are regular. If a singularity arises, it can point to a  breakdown of the  theory at short distances - as it happens for the equations of elasticity \cite{Gladwell}. Even if the theory does not break down,  singularities  might describe new fundamental phenomena, e.g., black holes or shock waves, which are unimaginable in the linear approximation. Thus, our understanding of a classical field theory is not complete until we know the conditions on the initial configuration under which regular solutions exist.

It is not usually difficult to show the existence of weak solutions, i.e., distribution solutions for which the equations hold  when averaged by a smooth test function. However, to prove that solutions are smooth  (``classical"\footnote{In the analysis of PDEs as well as in physics, the word ``classical" has a specific meaning; but the meaning is completely different in each discipline. To a physicist, a ``classical theory" is an approximation in which the quantum effects are ignored; classical equations of motion are often PDEs. In analysis, a ``classical solution" of a PDE is one that can  be represented as a function, not just a ``weak solution" which is only a distribution.})  is often  a very challenging mathematical problem. Moreover, to show that an evolution equation has a local solution, for some finite time after the initial condition, is often much easier than to show  global regularity.

Even for linear equations, the understanding of the regularity of solutions came many decades after physicists had already moved on to other matters. For example, the proof of regularity for solutions of the Laplace equation needed deep developments in functional analysis, such as Sobolev spaces, in the 1930s.

The most celebrated open problem of regularity is the one for the Navier-Stokes equations of hydrodynamics. While this ``millenium problem" has resisted many attempts, there has been a great deal of recent  progress in understanding the regularity of solutions for nonlinear wave equations, many of which also appear in physics. The simplest kind of nonlinear wave equation  is
\begin{equation}
g^{\mu\nu}\,\partial_\mu\partial_\nu \phi^i \,+\, {\partial V(\phi)\over \partial \phi^i}\,=\,0
\end{equation}
where the metric $g$ has signature $(n,1)$ and $V$ is some polynomial. This equation follows from the variational principle with action
\begin{equation}
S \,=\,  \int  \left[ {1\over 2} g^{\mu\nu} \partial_\mu\phi^i \partial_\nu\phi^i \, + \, V(\phi)\right]\,dx
\end{equation}
where the field $\phi:\mathbb{R}^{n,1}\to \mathbb{R}^m$ takes values into an Euclidean space .

Regular solutions for smooth initial data exist \cite{MR1078267,MR1162370} when $V$ is bounded below and its degree is not too large,
\begin{equation}
{\rm deg}\,V\,\leq\, {2n\over n-2}
\end{equation}
The inequality is saturated when the energy is scale invariant.  Thus, $V(\phi)\sim \phi^4$ is regular for spaces of dimension $n\leq 4$. When the dimension is larger, the theory is said to be supercritical and the solution is expected to develop singularities even for smooth initial data. For dimensions which are less than the critical one, the proof of regularity is substantially simpler.

It is interesting to note that a $\phi^4$  scalar quantum field theory is perturbatively renormalizable precisely when $n=3$ (the {\it spacetime} is four dimensional); i.e.,  one dimension lower than the condition for the classical theory to be regular. This is the condition for the action (not the energy) to be scale-invariant. Again, when the dimension is smaller than the critical one, the proof that the divergences of the Feynman diagrams can be removed is substantially simpler: there are no primitively divergent diagrams and the theory is super-renormalizable. It is only in this case that we have a mathematically rigorous construction of quantum field theories \cite{GlimmJaffe}. Thus, there appears to be some deep connection between the renormalizability of quantum field theories of some given space-time dimension and the regularity of the solutions for its classical equations of motion in one dimension higher.

An example of a renormalizable theory is the linear $\sigma$ model \cite{BenjaminLee} with four components $(\sigma,\pi_1,\pi_2,\pi_3)$ and
\begin{equation}
V(\sigma,\pi)={\lambda \over 4}(\sigma^2+\pi_1^2+\pi_2^2+\pi_3^2-a^2)^2+\epsilon\sigma
\end{equation}
The fields $\sigma$ and $\pi$ represent certain particles called mesons. The $\sigma$ meson is much heavier than the $\pi$ mesons, so that an approximation in which the pions are massless and the $\sigma$ is infinitely heavy is reasonable. As it happens, this is the limit of $\epsilon\to 0$ and $\lambda\to \infty$, keeping $a$ fixed. It turns out that this limiting process is in fact a better description of nature than the linear $\sigma$ model. We observe $\pi$ mesons everyday, whereas $\sigma$ mesons, if they exist at all, are very short-lived resonances.

In this limit we get the ``nonlinear $\sigma$ model" for which the field is constrained to lie on a three-sphere of radius $a$. If we introduce some coordinates on the sphere, the limiting equations can be written as
\begin{equation}
S\,=\,\frac 12 \int  g^{\mu\nu} \,\partial_\mu\phi^i\, \partial_\nu\phi^j\,
h_{ij}(\phi)\,dg
\end{equation}
where $h_{ij}$ is the metric of the sphere and $dg=\sqrt{-g}\, dx$ is the usual invariant measure in Minkowski spacetime. The Euler-Lagrange equation of this action is called the ``wave map equation" in the mathematics literature: it is a generalization of the geodesic equation to several indepedent variables, as well as the analogue of the harmonic map equation for a domain with Lorentzian signature.

The nonlinear $\sigma$ model is not renormalizable in four space-time dimensions (the critical dimension of its spacetime is two). Interestingly enough, regularity for the solutions of the wave map equation was shown in 2+1 dimensions \cite{G, SZ1, SZ2}. This is in keeping with our earlier observation that  the critical dimension of a quantum theory is one lower than that of the classical theory. In particular, we should expect that the wave map equation in four space-time dimensions has singular solutions even for regular initial data. This has been believed by physicists for a long time, based on scaling arguments. Specifically, the initial configurations that contain a ``topological soliton" evolve to a singularity in finite time. Since our paper concerns a possible solution to this pathology, let us digress to explain this point further.

To have finite energy, any initial data for the wave map equation must tend to a constant at spatial infinity. Identifying the points at infinity, this defines a continuous map from $\mathbb{S}^3$ to $\mathbb{S}^3$. An initial configuration of topological degree not equal to zero cannot be continuously deformed to a constant: it represents a new kind of excitation for the theory, fundamentally different from $\pi$ mesons, which are just infinitesimal perturbations. For small times, such a topologically nontrivial configuration will move like a particle of finite size. It is of much interest to consruct a theory in which such topological solitons exist and are stable (i.e., free of singularities). Alas, the nonlinear $\sigma$ model (the wave map equation) is not that theory: scaling arguments show that it is possible to reduce the energy of a topologically nontrivial configuration by shrinking it. As the size goes to zero, the energy goes to zero. Thus, an initial configuration of degree one should evolve in finite time to a singularity, where the degree changes discontinuously to zero. This is by now well established, indeed an explicit solution with such  a singularity can be exhibited \cite{TS}:
 \begin{equation}
\phi(t,x)\,=\,\left(\frac{2tx}{t^2+|x|^2}\,, \frac{|x^2|-t^2}{t^2+|x|^2}\right)
\end{equation}

But we should not give up on the wave map equation in three space dimensions entirely. The same physics that predicts the collapse of the soliton also suggests a way out: change the theory so that there is a short-distance repulsion which will halt the shrinking of the soliton at some finite size. Skyrme \cite{MR0128862} proposed a model that has exactly this feature. He further made the revolutionary suggesstion that the topological solitons of the nonlinear $\sigma$ model are atomic nuclei: the degree of the map is precisely the total number of protons and neutrons, a conserved quantity. These strange ideas have been found to be valid as an effective model \cite{SyracuseSkyrme,Witten} of the fundamental theory of nuclear interactions, Quantum Chromo-Dynamics (QCD). The profound connections between the Skyrme model and the topological anomalies of QCD, discovered by Witten \cite{Witten}, were crucial in establishing this relation. An exact equivalence  to QCD can be proved in  two space-time dimensions \cite{RajeevQHD}.

Physics suggests that the wave map equation modified by Skyrme's repulsion will have regular solutions even for initial data that are of non-zero degree. However, proving this is a difficult problem: the Skyrme equation is even more nonlinear than the wave map equation. More precisely, it is no longer semilinear, but is  quasilinear. This means that the coefficients of the top order terms are functions of the field, and, as a result, not even the local existence theory follows from standard theorems, but has to be proven directly.

There exists however a simpler  alternative to the Skyrme repulsion, suggested by Adkins and Nappi \cite{AN}, which leads to a semilinear equation. Physically, it describes the interactions of three $\pi$ mesons with an $\omega$ vector meson field, which provides the repulsive force. Scaling arguments, as well as the existence of  a static solution, suggest that this equation has regular solutions even for large initial data. We thus conjecture that the Adkins-Nappi model is no longer supercritical, i.e., solutions starting from finite energy smooth initial data are global in time\footnote{The limit of zero mass for the $\omega$ and $\pi$ mesons is much simpler to study.  It is expected that the short-distance behavior, such as singularities, are unaffected by the masses of the mesons. This is why we will study the massless limit of the Adkins-Nappi model.}.

By now there is a well-established approach to proving such regularity theorems, which has worked for the wave map equation in lower dimensions. Due to the fact that our equation is semilinear, we know right away that for any smooth data there exists at least a local solution.  First, one proves that energy does not concentrate (i.e., the amount of energy within a ball tends to zero as its  radius goes to zero) and that the solution remain continuous at the first possible singularity. If we can also show that small energy implies regularity, we would have achieved our goal of proving regularity.

Earlier \cite{GR1} we succeded in showing that a spherically symmetric (equivariant) solution is continuous. In this paper, we complete the first part of the above scheme by showing nonconcentration of the energy, again, in the case of a spherically symmetric  solution. As time evolves, some energy is radiated away and what remains should settle into a static solution.

\section{Main results}

The model we are studying is given by the action
\begin{equation}
S\,=\,\frac 12 \int  g^{\mu\nu} \,\partial_\mu\phi^i\, \partial_\nu\phi^j\,
h_{ij}(\phi)\,dg\,+\, \frac{1}{4} \int F^{\mu\nu} F_{\mu\nu}\,dg\,-\, \int A_\mu j^\mu \, dx
\label{S}
\end{equation}
where:
\begin{itemize}

\item \,$\phi :(\mathbb{R}^{3+1},g)\to(\mathbb{S}^3,h)$ is a map from
the Minkowski spacetime, with $g= \text{diag}(-1,1,1,1)$, into
the sphere $\mathbb{S}^3\subset\mathbb{R}^{4} $ (with the
induced Riemannian metric $h$), describing the $\pi$ mesons;

\item \,the 1-form $A=A_\mu dx^\mu$ is a gauge potential representing the $\omega$ meson, while the 2-form
$F_{\mu\nu}\,=\,\partial_\mu A_\nu -
\partial_\nu A_\mu$ is its associated electromagnetic field;

\item \,$j$ is the flux or the baryonic current, i.e.,
\[
j^\mu\,=\,c\, \epsilon^{\mu\nu\rho\sigma}\,
\partial_\nu\phi^i \,\partial_\rho\phi^j \,\partial_\sigma\phi^k
\,\epsilon_{ijk}
\]
with $\epsilon$ being the Levi-Civita symbol and $c$ a normalizing constant.
\end{itemize}

The first term in \eqref{S} is the action of the nonlinear $\sigma$ model, while the third one describes the coupling of the $\pi$ mesons with the vector meson.

For a field $\phi$ of finite energy, one can associate naturally its topological charge (winding number)
\begin{equation}
Q=\int c\, \epsilon^{\nu\rho\sigma}\,
\partial_\nu\phi^i \,\partial_\rho\phi^j \,\partial_\sigma\phi^k
\,\epsilon_{ijk}\,dx\,=\,\int j^0 \,dx \label{Q}
\end{equation}
In the Adkins-Nappi model, the density of the gauge potential $A$ coincides with the density of the topological charge $Q$.

The object of our investigation are time dependent equivariant maps of winding number 1 associated to \eqref{S}, i.e.,
\[
\phi(t,r,\psi, \theta)= (u(t,r),\psi, \theta)\qquad A(t,r,\psi,
\theta)=(V(t,r),0,0,0)
\]
\[
u(t,0)=0  \qquad u(t,\infty)= \pi
\]
that yield the following PDE for $u$:
\begin{equation}
u_{tt} - u_{rr} - \frac 2r u_r + \frac{\sin 2u}{r^2} +
\frac{(u-\sin u \cos u)(1- \cos 2u)}{r^4} = 0\label{main}
\end{equation}

A simple calculus shows that the energy norm
\begin{equation}
\mathcal{E}[u](t)=\int_0^\infty \left[\frac{1}{2}\left(u_t^2+u_r^2\right)+
\frac{\sin ^2 u}{r^2} +\frac{(u-\sin u \cos u)^2}{2r^4}
\right]\,r^2 dr \label{tote}
\end{equation}
is a conserved quantity if $u$ satisfies \eqref{main}.

Based on the fact that  \eqref{main} is a semilinear equation, for which smooth solutions exist at least locally, and $u$ is a radial function, we can assume, without any loss of generality, that our solution starts at time $T_0= -1$ and develops its first singularity at $(0,0)$. In this
context, our main result is

\begin{theorem}
The local energy of a solution for \eqref{main} does not concentrate:
\begin{equation}
\lim_{T\to 0-} \int_0^{|
T|}\,\left[\frac{1}{2}\left(u_t^2+u_r^2\right)+
\frac{\sin ^2 u}{r^2} +\frac{(u-\sin u \cos u)^2}{2r^4}
\right]\, r^2\,dr\,=\,0\label{en}
\end{equation}
\label{local}
\end{theorem}

This should be compared with what we were able to prove in \cite{GR1}:

\begin{theorem}
(\cite{GR1}) If $u$ satisfies \eqref{main}, then the following local energy estimate holds
\begin{equation}
\lim_{T\to 0-} \int_0^{|
T|}\left[ \left(1-\frac{r}{|T|}\right)u_-^2+u_+^2 +\frac{\sin^2u}{r^2} +\frac{(u-\sin
u \cos u)^2}{r^4}\right] r^2 dr = 0\label{loce}
\end{equation}
where
\[
u_{+}=u_t+u_r \qquad u_{-}=u_t-u_r
\]
is the standard notation for derivatives along null directions. As a consequence, $u$ is continuous at the origin with
\[
\lim_{(t,r)\to (0,0)} u(t,r)\,=\, 0
\]
\label{cont}
\end{theorem}

\begin{remark}
The continuity argument uses only
\[
\lim_{T\to 0-} \int_0^{|
T|}\,\frac{(u-\sin
u \cos u)^2}{r^4}\, r^2 \,dr \,=\, 0\]
\end{remark}

\begin{remark}
One notices that what we are missing in \eqref{loce} from claiming \eqref{en} is an estimate near the cone (i.e., the region where $r \approx |T|$) for the null direction derivative $u_-$.
\end{remark}

In proving Theorem \ref{cont}, we used the conventional strategy (e.g., \cite{SZ1}, \cite{MR1223662}) of obtaining local energy estimates by integrating differential identities  over appropriate past oriented domains. These were derived in turn by multiplying \eqref{main} with various conformal Killing vector fields. However, this method
does not get us quite to where we want (i.e., proving \eqref{en}), mainly because, in our case, the energy is not critical (scale-invariant), some of its terms being supercritical\footnote{For a more detailed discussion of this aspect, we refer the reader to Remark 2.6 in \cite{GR1}}.

This is why here we adapt a technique pioneered by Grillakis \cite{G}, in which we work backwards in time (i.e., from positive times towards $0$) and take into account the behavior of $u$, both on past and future oriented domains. We would still be integrating differential identities, but Grillakis's method has the extra advantage of working with two apparently different expressions for the local energy:
\begin{equation}\aligned
E[u](T)\,&=\, \int_0^{
T}\left[\frac{1}{2}(u_t^2+u_r^2)+\frac{\sin^2u}{r^2}+\frac{(u-\sin
u \cos u)^2}{2r^4}\right] r^2\,dr\\
\,&=\, \frac{1}{\sqrt{2}}\int_{
K_T}\left[\frac{1}{2}(u_t-u_r)^2+\frac{\sin^2u}{r^2}+\frac{(u-\sin
u \cos u)^2}{2r^4}\right]
\endaligned
\label{equiv}
\end{equation}
where $K_T=\{(t,x)|\, T\leq t = 2T-|x|\leq 2T\}$. The first one is an integral over a spacelike ball, while the other is based on a light cone. They are equal because their domains form the boundary of a solid cone in spacetime. However, this equivalence gives better control on the behavior of energy on small balls.

\section{Notations, multipliers, and basic estimates}
We use the framework and notation of \cite{G}, in which $u$ is assumed to be a smooth solution for \eqref{main} in the domain
\begin{equation}
D(\overline{T})\,=\,\{(t,x)|\, |x|\leq\min\{2\overline{T}-t,t\},\,
0< t\leq 2\overline{T}\} \label{dom}
\end{equation}
with initial data given by
\begin{equation}
u(2\overline{T}-r,r)\,=\,f(r), \quad \overline{T}\leq r \leq
2\overline{T} \label{inidata}
\end{equation}

For $0<t_0<t_1\leq \overline{T}$, we denote
\[
\aligned D(t_0,t_1)\,&=\,\{(t,x)|\, 2t_0-t \leq
|x| \leq \min\{2t_1-t,t\}\}\\
C(t_0,t_1)\,&=\,\{(t,x)|\, t_0 \leq
t=|x| \leq t_1\}\\
\tilde{C}(t_0,t_1)\,&=\,\{(t,x)|\, |x|=t-2t_0 \leq
t_1-t_0\}\\
D_{t_0}\,&=\,\{(t_0,x)|\,
|x| \leq t_0\}\\
K_{t_0}\,&=\,\{(t,x)|\,  t_0 \leq t=2t_0-|x| \leq 2t_0\}
\endaligned
\]

\begin{center}
\includegraphics[width=3.5in]{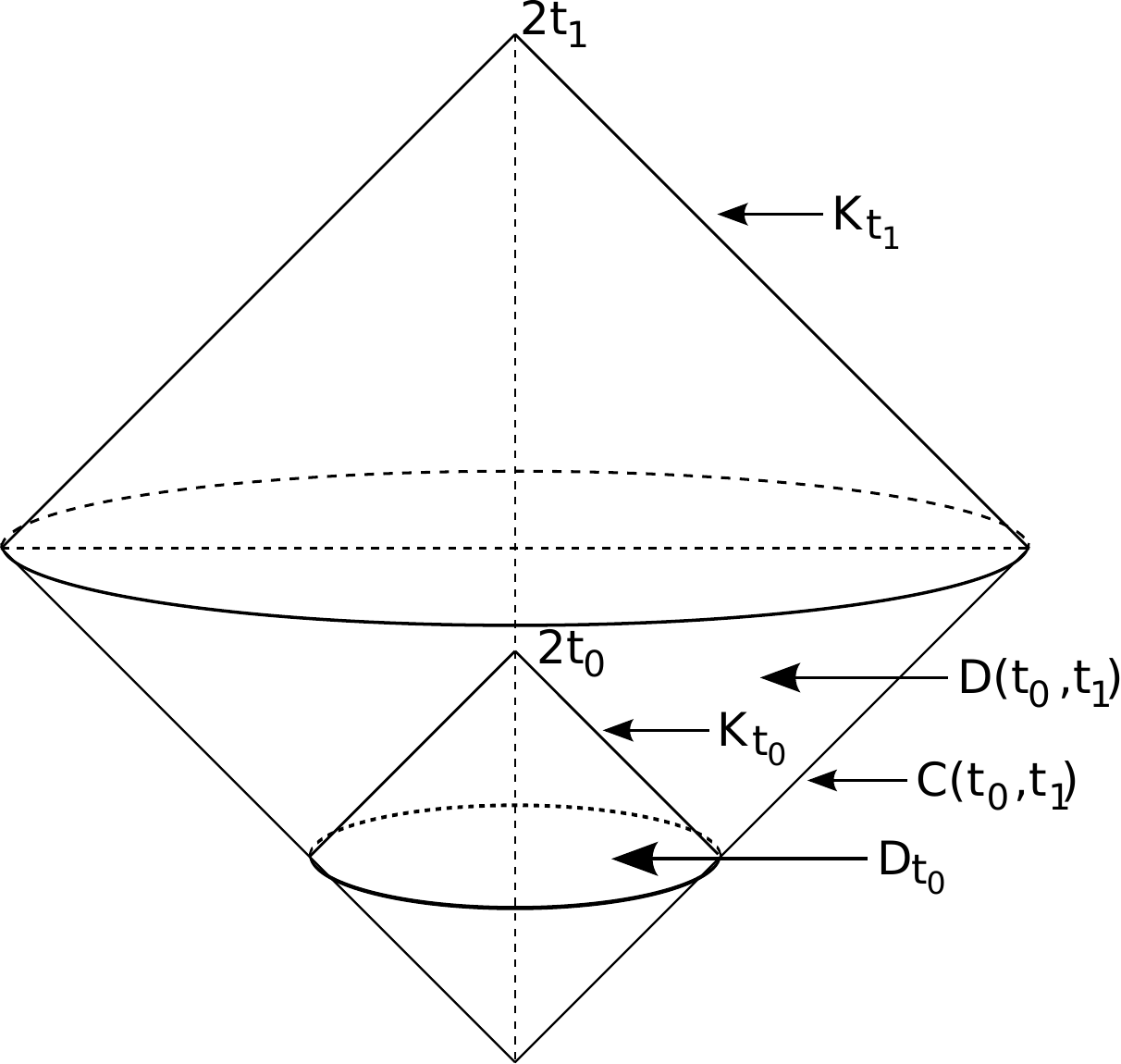}
\end{center}

Using the multiplier $a u_t + b u_r + c u$, where
$a$, $b$, and $c$ are all functions of $t$ and $r$, we obtain that $u$ satisfies:
\begin{equation}
\aligned &\partial_t\left( \frac{a+b}{4}\, u_{+}^2+
\frac{a-b}{4}\, u_{-}^2 + a\left(\frac{\sin^2u}{r^2}+\frac{(u-\sin
u \cos u)^2}{2r^4}\right) + cu u_t
-c_t\frac{u^2}{2}\right)\\
&-\frac{1}{r^2}\partial_r \bigg(r^2 \bigg[\frac{a+b}{4} u_{+}^2- \frac{a-b}{4}
u_{-}^2 - b\left(\frac{\sin^2u}{r^2}+\frac{(u-\sin u \cos
u)^2}{2r^4}\right) + cu u_r\\
&-c_r\frac{u^2}{2}\bigg]\bigg) = \left(\frac{a_t}{2}-\frac{\partial_r(br^2)}{2r^2} +
c\right) u_{t}^2+ (-a_r+b_t) u_t u_r +
\left(\frac{a_t-b_r}{2}+\frac{b}{r} - c\right) u_{r}^2 \\
& +(a_t + b_r)\frac{\sin^2u}{r^2}+\left(a_t +
r^2\partial_r\left(\frac{b}{r^2}\right)\right)\frac{(u-\sin u \cos
u)^2}{2r^4} - \frac{\Box c}{2}u^2 + c u \Box u= I
\endaligned \label{abc}
\end{equation}

For $(a,b,c)=(1,0,0)$ we derive the energy differential identity
\begin{equation}
\partial_t\left(\frac{1}{2}(u_t^2+u_r^2)+\frac{\sin^2u}{r^2}+\frac{(u-\sin
u \cos u)^2}{2r^4}\right)-\frac{1}{r^2}\,\partial_r(r^2u_tu_r)=0 \label{e1}
\end{equation}
which, integrated over the interior of the cone bounded by $K_t$ and $D_t$, justifies \eqref{equiv}.

Second, integrating \eqref{e1} on the domain $D(t_0,t_1)$, we obtain
\begin{equation}
E(t_1)-E(t_0)\,=\,F(t_0,t_1) \label{e}
\end{equation}
where $F(t_0,t_1)$ is the flux between the time slices $D_{t_0}$ and $D_{t_1}$, i.e.,
\[
F(t_0,t_1)\,=\,\frac{1}{\sqrt 2}\int_{C(t_0,t_1)}
\left[\frac{1}{2}\,u_+^2+\frac{\sin^2u}{r^2}+\frac{(u-\sin u \cos
u)^2}{2r^4}\right]
\]

Arguing exactly as in Section 2 of \cite{GR1}, where the problem was forward in time, we can infer from \eqref{e}:

\begin{prop}
In the framework detailed above (i.e., $u$ is a smooth solution for \eqref{main} in the domain $D(\overline{T})$), the following estimates hold:
\begin{align}
&E(t_0)\,\leq\, E(t_1), \qquad (\forall) 0< t_0\leq t_1 \leq
\overline{T} \label{mone}\\
&\lim_{T\to 0+}F(0,T)\,=\,\lim_{T\to 0+} u(T,T)\,=\,0
\label{flux+f0}\\
&\|u\|_{L^\infty}\leq C(E(\overline{T})) \label{li}
\end{align}
where $F(0,T)=\lim_{S\to 0} F(S,T)$ and $C(E(\overline{T})) $ is a constant depending strictly on $E(\overline{T})$.
\label{old}\end{prop}

Finally, integrating \eqref{e1} on the domain bounded by $K_{t_0}$, $C(t_0,t_1)$, $K_{t_1}$, and $\tilde{C}(t_0,t_1)$, and using \eqref{mone}-\eqref{flux+f0}, we deduce that
\[
\lim_{t_0\to 0} E(t_0) \,=\,\lim_{t_0\to 0} \frac{1}{\sqrt 2}\int_{\tilde{C}(t_0,t_1)}
\left[\frac{1}{2}\,u_+^2+\frac{\sin^2u}{r^2}+\frac{(u-\sin u \cos
u)^2}{2r^4}\right]\,-\,F(0,t_1)\]
which allows us to reduce the proof of \eqref{en} to the one of
\begin{equation}
\lim_{t_0\to 0} \frac{1}{\sqrt 2}\int_{\tilde{C}(t_0,t_1)}
\left[\frac{1}{2}\,u_+^2+\frac{\sin^2u}{r^2}+\frac{(u-\sin u \cos
u)^2}{2r^4}\right]\,=\,F(0,t_1)\label{lim}
\end{equation}
for $t_1$ sufficiently small.

\section{Local energy arguments}
We proceed by integrating \eqref{abc} on the domain $D(t_0,t_1)$ to obtain:
\begin{equation}
\aligned &\sqrt 2\,\int_{D(t_0,t_1)} I + \int_{C(t_0,t_1)} \left[ G(u_+, \sin u, u-\sin u \cos u)-\left(c_r+c_t+\frac{c}{r}\right)u^2\right]\\
&=\,\int_{K_{t_1}} - \int_{K_{t_0}} \left[ H(u_-, \sin u, u-\sin u \cos u)+\left(c_r-c_t+\frac{c}{r}\right)u^2\right]\\
\endaligned\label{int}
\end{equation}
where
\[
\aligned
G(u_+, \sin u, u-\sin u \cos u)\,&=\,\frac{a+b}{2}\, u_{+}^2+
(a-b)\left(\frac{\sin^2u}{r^2}+\frac{(u-\sin u \cos
u)^2}{2r^4}\right)\\
H(u_-, \sin u, u-\sin u \cos u)\,&=\,\frac{a-b}{2}\, u_{-}^2+
(a+b)\left(\frac{\sin^2u}{r^2}+\frac{(u-\sin u \cos
u)^2}{2r^4}\right)\endaligned
\]

To go beyond the basic estimates \eqref{mone}-\eqref{li} and obtain sharper ones, that would allow us to derive \eqref{lim}, the first challenge is to find multipliers, other than $(a,b,c)=(1,0,0)$, for which the integrals on $K_{t_1}$, $K_{t_0}$, and $C(t_0,t_1)$, are comparable with the energy, respectively the flux. As we control the energy, the flux, and the size of $u$ through \eqref{mone}-\eqref{li}, the second goal would be to show that $\int_{D(t_0,t_1)} I $ is uniformly bounded for $0<t_0<t_1\leq \overline{T}$.

We can achieve these by using $(a,b,c)\,=\,\left(1, \frac{r}{t}, \frac{1}{t} \right)$, for which one deduces\footnote{In what follows, the notation $A\lesssim B$ will stand for $A\leq C B$, where $C$ is a constant independent of $A$ and $B$.}
\begin{equation}
\aligned
0\leq G(u_+, \sin u, u-\sin u \cos u) &\lesssim \frac{1}{2}\,u_+^2+\frac{\sin^2u}{r^2}+\frac{(u-\sin u \cos
u)^2}{2r^4}\\
0\leq H(u_-, \sin u, u-\sin u \cos u) &\lesssim \frac{1}{2}\,u_-^2+\frac{\sin^2u}{r^2}+\frac{(u-\sin
u \cos u)^2}{2r^4}\\
\left| c_r+c_t+\frac{c}{r} \right| \lesssim \frac{1}{r^2}, &\qquad 0 \leq c_r-c_t+\frac{c}{r} \lesssim \frac{1}{r^2}
\endaligned
\label{gh}
\end{equation}
Taking into account also  \eqref{e}, \eqref{mone}, and \eqref{li}, we are lead to
\begin{equation}
\left|\int_{D(t_0,t_1)} I \right| \leq  C(E(\overline{T})), \qquad 0<t_0<t_1\leq \overline{T}\label{bi}
\end{equation}

Next, we want to prove that the above estimate implies energy decay. We show initially that one can ignore a number of terms from $\int_{D(t_0,t_1)} I $, which are also uniformly bounded. For $(a,b,c)\,=\,\left(1, \frac{r}{t}, \frac{1}{t} \right)$:
\begin{equation}
\left| \int_{D(t_0,t_1)} \left[(a_t + b_r)\frac{\sin^2u}{r^2} - \frac{\Box c}{2}u^2 - c\,\frac{u\sin 2u}{r^2}\right] \right| \leq  C(\overline{T}, E(\overline{T})) ,\, 0<t_0<t_1\leq \overline{T}
\label{bd}
\end{equation}
which follows from \eqref{li} and
\[
\int_{D(t_0,t_1)} \frac{1}{t\,r^2} \approx \int_0^{t_0}\int_{2t_0-r}^{2t_1-r} + \int_{t_0}^{t_1}\int_{r}^{2t_1-r} \, \frac{1}{t}\,dt\,dr\,\lesssim \, \overline{T} + \sup_{[0,\overline{T}]} |x \log x|
\]

Coupling \eqref{bi}-\eqref{bd} with
\[\aligned
\left(\frac{a_t}{2}-\frac{\partial_r(br^2)}{2r^2} +
c\right) u_{t}^2+\, &(-a_r+b_t) u_t u_r +
\left(\frac{a_t-b_r}{2}+\frac{b}{r} - c\right) u_{r}^2
\lesssim \frac{b-a}{t} u_-^2\\
&a_t +
r^2\partial_r\left(\frac{b}{r^2}\right)\,=\,-\frac{1}{t}
\endaligned\]
which are obtained through direct computations, and the crucial fact\footnote{used also in \cite{GR1}}
\begin{equation}
u(u-\sin u \cos u)(1- \cos 2u)\,\geq 0\label{pos}\end{equation}
we derive:
\[
\int_{D(t_0,t_1)} \frac{1}{t}\left[(1-\frac{r}{t})u_-^2 + \frac{(u-\sin u \cos u)^2}{r^4}\right] \leq  C(\overline{T}, E(\overline{T})), \ \ 0<t_0<t_1\leq \overline{T}
\]
This yields the existence of a sequence $(s_n)_n \to 0$ for which
\begin{equation}
\lim_{n\to\infty} \int_{K_{s_n}}  (1-\frac{r}{t})u_-^2 + \frac{(u-\sin u \cos u)^2}{r^4}\,=\,0
\label{de}
\end{equation}

Hence, we are allowed to take $t_0 = 0$ in \eqref{int} (for $(a,b,c)\,=\,\left(1, \frac{r}{t}, \frac{1}{t} \right)$) to obtain:
\begin{equation}
\aligned &\sqrt 2\,\int_{D(0,T)} I + \int_{C(0,T)} \left[ G(u_+, \sin u, u-\sin u \cos u)-\left(c_r+c_t+\frac{c}{r}\right)u^2\right]\\
&=\,\int_{K_{T}}  \left[ H(u_-, \sin u, u-\sin u \cos u)+\left(c_r-c_t+\frac{c}{r}\right)u^2\right], \quad 0<T\leq \overline{T}\\
\endaligned
\label{fint}
\end{equation}
In deducing \eqref{fint} we relied also on the trivial observation
\begin{equation}
\lim_{t\to 0}  \int_{C(0,t)} \frac{f}{r^2} \,=\,
\lim_{t\to 0} \int_{K_t} \frac{f}{r^2} \,=\,0 \quad \text{for} \ \ f \in L^\infty
\label{r2}
\end{equation}

Using the decay of the flux \eqref{flux+f0} for the integral on $C(0,T)$ and a slight variation of \eqref{bd}, i.e.,
\[
\lim_{T\to 0} \int_{D(0,T)} (a_t + b_r)\frac{\sin^2u}{r^2} - \frac{\Box c}{2}u^2 - c\,\frac{u\sin 2u}{r^2} \,=\,0
\]
we can infer from \eqref{fint} that

\begin{prop}
Under the hypotheses of Proposition \ref{old}, $u$ satisfies the energy decay
\begin{equation}
\lim_{T\to 0}\, \int_{K_{T}} \left(1-\frac{r}{t}\right) u_-^2 + \frac{\sin^2u}{r^2}+\frac{(u-\sin u \cos u)^2}{r^4}\,=\,0
\label{energy}\end{equation}
\end{prop}

\begin{remark}
One needs to compare \eqref{energy} with  \eqref{loce} (obtained in \cite{GR1}) and realize that \eqref{energy} is also enough to clinch the continuity argument.
\end{remark}

We use now \eqref{int} for $(a,b,c)=(1,1,\frac 1r)$ to derive
\[
\aligned
-\sqrt 2\,\int_{D(t_0,t_1)}  &\frac{(u-\sin u \cos u)^2}{r^5} + \frac{u\sin{2u}}{r^3} + \frac{u(u-\sin u \cos u)(1- \cos 2u)}{r^5}\\
&=\,\int_{K_{t_1}} - \int_{K_{t_0}} \frac{2\sin^2u}{r^2}+\frac{(u-\sin u \cos u)^2}{r^4} - \int_{C(t_0,t_1)} u_+^2\endaligned\]
Relying on \eqref{pos}, \eqref{energy}, and that $u\sin{2u}\geq 0$ near the origin, one deduces
\begin{equation}
\int_{D(0,T)} \frac{(u-\sin u \cos u)^2}{r^5} + \frac{u\sin{2u}}{r^3} + \frac{u(u-\sin u \cos u)(1- \cos 2u)}{r^5} < +\infty
\label{d0t}\end{equation}
for $T$ sufficiently small.

We can then integrate \eqref{abc}, for $(a,b,c)=(1,1,\frac 1r)$, on the domain $\tilde{D}(t_0,t_1)$, bounded by $\tilde{C}(t_0,t_1)$, $K_{t_1}$, and $C(0,t_1)$
\begin{center}
\includegraphics[width=3.5in]{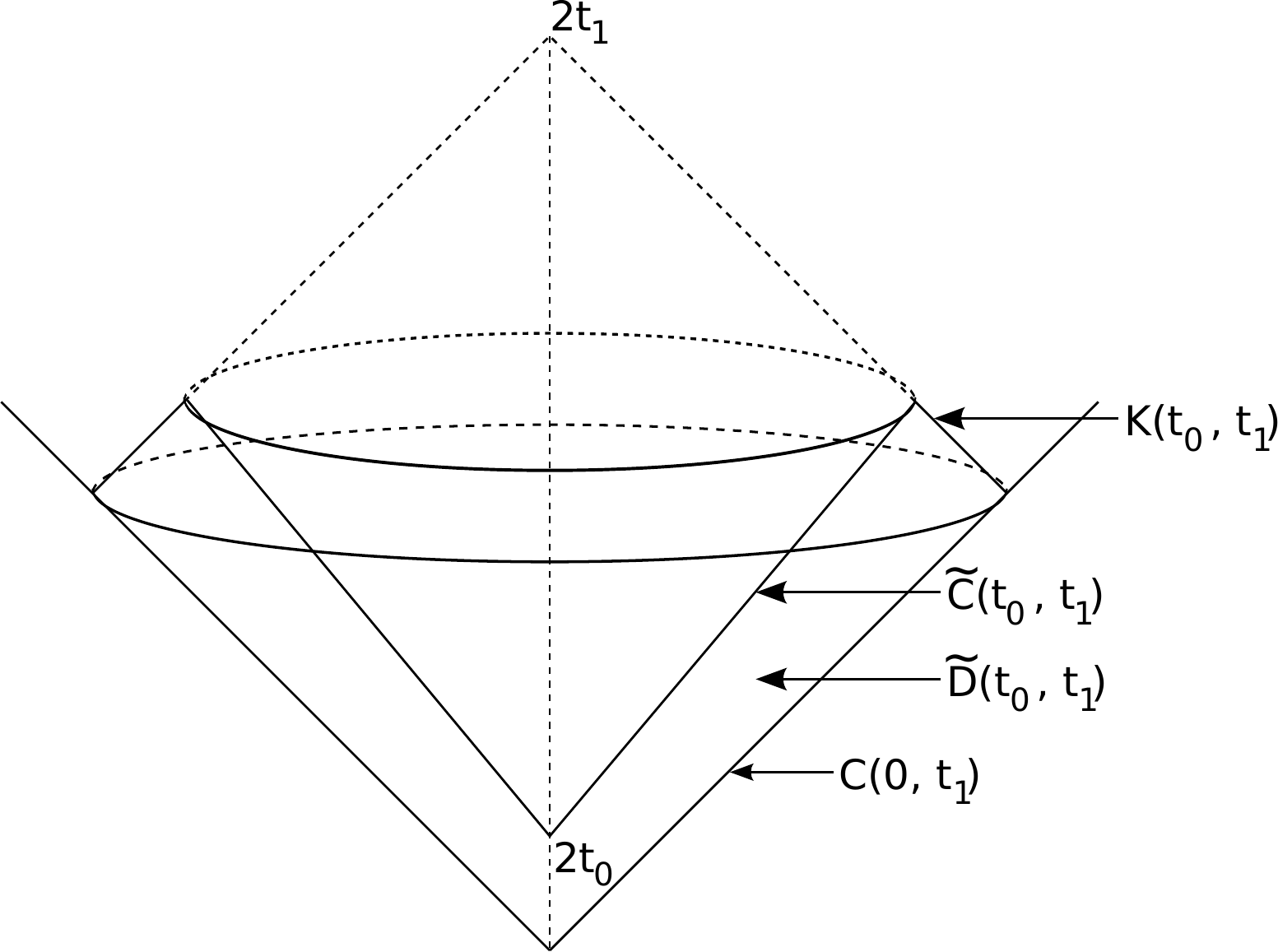}
\end{center}
to obtain
\[
\aligned
&\int_{C(0,t_1)} - \int_{\tilde{C}(t_0,t_1)} u_+^2\,=\,\int_{K(t_0,t_1)} \frac{2\sin^2u}{r^2}+\frac{(u-\sin u \cos u)^2}{r^4} \\
+ \,\sqrt 2\,&\int_{\tilde{D}(t_0,t_1)}  \frac{(u-\sin u \cos u)^2}{r^5} + \frac{u\sin{2u}}{r^3} + \frac{u(u-\sin u \cos u)(1- \cos 2u)}{r^5}
\endaligned
\]
where
\[
K(t_0,t_1)\,=\,\{(t,x)|\,  t_1 \leq t=2t_1-|x| \leq t_1+t_0\}\,\subseteq \,K(t_1)
\]
The integrals on $K(t_0,t_1)$ and $\tilde{D}(t_0,t_1)$ are handled through \eqref{mone} and \eqref{d0t}, and so we can imply
\[\lim_{t_0\to 0} \int_{\tilde{C}(t_0,t_1)} u_+^2\,=\,\int_{C(0,t_1)} u_+^2\]
for $t_1$ sufficiently small.

A similar argument for $(a,b,c)=(1,\frac 12,\frac{1}{2r})$ leads to
\[
\aligned
\lim_{t_0\to 0} \int_{\tilde{C}(t_0,t_1)} &\frac{3}{4}u_+^2 +  \frac{\sin^2u}{2r^2}+\frac{(u-\sin u \cos u)^2}{4r^4}=\\
&= \int_{C(0,t_1)} \frac{3}{4}u_+^2 +  \frac{\sin^2u}{2r^2}+\frac{(u-\sin u \cos u)^2}{4r^4}
\endaligned
\]
Combining these last two limits, we derive \eqref{lim} and thus finish the argument.

\section*{Acknowledgements}
We thank Manoussos Grillakis for stimulating discussions in the early stages of this work and our colleague, Scott Bailey, for helping us with the figures. The first author was supported in part by the National Science Foundation Career grant DMS-0747656. The second author was supported in part by the Department of Energy contract DE-FG02-91ER40685.

\bibliographystyle{amsplain}

\end{document}